\documentclass[11pt]{amsart}
\usepackage{amssymb}
\usepackage{amsmath}

\numberwithin{equation}{section}

\newtheorem{thm}{Theorem}[section]
\newtheorem{lem}[thm]{Lemma}
\newtheorem{cor}[thm]{Corollary}

\newcommand{\be}{\begin{equation}}
\newcommand{\ee}{\end{equation}}
\newcommand{\bes}{\begin{equation*}}
\newcommand{\ees}{\end{equation*}}

\newcommand{\bea}{\begin{eqnarray}}
\newcommand{\eea}{\end{eqnarray}}
\newcommand{\ba}{\begin{array}}
\newcommand{\ea}{\end{array}}
\newcommand{\bc}{\begin{center}}
\newcommand{\ec}{\end{center}}

\def\s{\sigma}
\def\e{{\bf 1}\!\!{\rm I}}
\def\l{\lambda}

\def\i{\varepsilon}
\def\t{\tau}
\def\f{\varphi}
\def\N{\mathbb{N}}
\def\R{\mathbb{R}}
\def\a{\alpha}

\def\m{\mu}
\def\n{\nu}

\def\cf{{\mathcal F}}

\begin{document}
\title[ON MIXING AND COMPLETELY MIXING PROPERTIES]
{ON MIXING AND COMPLETELY MIXING PROPERTIES OF POSITIVE
$L^1$-CONTRACTIONS OF FINITE VON NEUMANN ALGEBRAS}
\thanks{The
work supported by NATO-TUBITAK PC-B programme}
\author{Farrukh Mukhamedov}
\address{Farrukh Mukhamedov\\
Department of Mechanics and Mathematics\\
National University of Uzbekistan\\
Vuzgorodok, 700174, Tashkent, Uzbekistan} \email{{\tt
far75m@yandex.ru}}
\author{Seyit Temir}
\address{Seyit Temir\\
Department of Mathematics\\
 Arts and Science Faculty\\
 Harran University, 63200 \\
 \c{S}anliurfa, Turkey}
\email{{\tt seyittemir67@hotmail.com}}
\author{Hasan Akin}
\address{Hasan Akin\\
Department of Mathematics\\
 Arts and Science Faculty\\
 Harran University, 63200 \\
 \c{S}anliurfa, Turkey}
\email{{\tt hasanakin69@hotmail.com}}



\begin{abstract}
 Akcoglu and Suchaston proved
the following result: Let $T:L^1(X,{\cf},\m)\to L^1(X,{\cf},\m)$
be a positive contraction. Assume that  for $z\in L^1(X,{\cf},\m)$
the sequence $(T^nz)$ converges weakly in $L^1(X,{\cf},\m)$, then
either $\lim\limits_{n\to\infty}\|T^nz\|=0$ or there exists a
positive function $h\in L^1(X,{\cf},\m)$, $h\neq 0$ such that
$Th=h$. In the paper we prove an extension of this result  in
finite von Neumann algebra setting, and as a consequence we obtain
that if a positive contraction of a noncommutative  $L^1$-space
has no non zero positive invariant element, then its mixing
property implies completely mixing property one.\\[2mm]
\end{abstract}

\maketitle


\section{Introduction}

It is known (see \cite{K}) that there are several notions of
mixing (i.e. weak mixing, mixing, completely mixing e.c.t.) of
measure preserving transformation on a measure space in the
ergodic theory. It is important to know how these notions are
related with each other. A lot of papers are devoted to this
topic. For example, recently, in \cite{BKLM} relations between the
notions of weak mixing and weak wandering have been studied.

In this paper we deal with the notions of mixing and completely
mixing. Now recall them. Let $(X,\cf,\m)$ be a measure space with
probability measure $\m$. Let $L^1(X,\cf,\m)$ be the associated
$L^1$-space. A linear operator $T:L^1(X,\cf,\m)\to L^1(X,\cf,\m)$
is called {\it positive contraction} if  $Tf\geq 0$ whenever
$f\geq 0$ and $\|T\|\leq 1$. Let $L^1_0=\{f\in L^1(X,\cf,\m): \int
fd\m=0\}$. A positive contraction $T$ in $L^1(X,{\cf},\m)$ is
called {\it mixing} (resp. {\it completely mixing}) if $T^nf$
tends weakly to 0 for all $f\in L^1(X,{\cf},\m)$ (resp. $\|T^nf\|$
tends to 0 for all $f\in L^1_0$). A relation between these two
notions was given in (see \cite{KS}), which can be formulated as
follows
(see, \cite{K}, Ch.8, Th-m 1.4):\\

\begin{thm}\label{ks} Let $T:L^1(X,{\cf},\m)\to
L^1(X,{\cf},\m)$ be a positive contraction. Assume that there
exists no non zero $y\in L^1(X,{\cf},\m)$, $y\geq 0$ such that
$Ty=y$. If for $z\in L^1(X,{\cf},\m)$ the sequence $(T^nz)$
converges weakly to some element of $L^1(X,{\cf},\m)$, then
$\lim\limits_{n\to\infty}\|T^nz\|=0$. In particular, if $T$ is
mixing, then $T$ is completely mixing. \end{thm}

This theorem gives an answer to the problem whether
K-automorphisms of $\s$-finite measure space are mixing, and
showed that, in fact, invertible mixing measure preserving
transformations of $\s$-finite infinite space do not exist (see
\cite{KS}) (see for review \cite{K}).

Later in \cite{AS} Akcoglu and Sucheston proved an extension of
Theorem \ref{ks}, which is formulated as follows

\begin{thm}\label{as}  Let $T:L^1(X,{\cf},\m)\to L^1(X,{\cf},\m)$ be a positive
contraction.  Assume that  for $z\in L^1(X,{\cf},\m)$ the sequence
$(T^nz)$ converges weakly in $L^1(X,{\cf},\m)$, then either
$\lim\limits_{n\to\infty}\|T^nz\|=0$ or there exists a positive
function $h\in L^1(X,{\cf},\m)$, $h\neq 0$ such that $Th=h$.
\end{thm}

By means of this theorem in \cite{AS},\cite{FS} an extension of
Blum-Hanson theorem was proved \cite{BH}, which states that if $T$
is a positive contraction on $L^1(X,{\cf},\m)$, then $T$ is mixing
if and only if $\frac{1}{n}\sum\limits_{i=1}^nT^{k_i}f$ converges
strongly in $L^1(X,{\cf},\m)$ for every $f\in L^1(X,{\cf},\m)$ and
increasing sequence of integers $\{k_n\}$: $0\leq k_1<k_2\cdots$.
Other extensions of this result have been given in \cite{BLRT}.

The formulated Theorem \ref{as} has a lot of applications, here we
mention only a few of them. Namely, using it in \cite{E} the
existence of an invariant measure for given positive contraction
$T$ on $L^1(X,{\cf},\m)$ was proved and in \cite{ZZ} a criterion
of strong asymptotically stability for positive contractions was
given by means of Theorem \ref{as}.

In this paper we are going to extend these results for quantum
dynamical systems over von Neumann algebras. Here by quantum
dynamical systems we mean a linear, positive, weak continuous
mapping $\a$ of a von Neumann algebra $M$, with normal faithful
trace $\t$, into itself. It is known (see \cite{BR}, sec.4.3) that
the theory of quantum dynamical systems provides convenient
mathematical description of the irreversible dynamics of an open
quantum system. This motivates an interest in the study of
conditions for a dynamical system to induce approach to a
stationary state,  of reflect subjects such as irreducibility
(i.e. ergodicity, mixing) and ergodic theorems (see for example,
\cite{AH},\cite{FR1},\cite{J1}). By means of given dynamical
system $\a$ a conjugate dynamical system $\a_*$ on a predual $M_*$
of the algebra $M$ can be defined  as follows
$(\a_*\f)(x)=\f(\a(x))$, where $\f\in M_*$,$x\in M$. It is known
\cite{N} that a predual $M_*$ of $M$ is isometrically isomorphic
to $L^1(M,\t)$, which is a non-commutative analog of $L^1$-space.
Note that the study of limiting behaviors of $\a_*$ is important
since, it will provide towards a proof of some ergodic type
theorems for $\a$. Therefore, in the paper we consider only
$L^1$-contractions of $L^1(M,\t)$. Note that such kind of
dynamical systems were considered in \cite{S},\cite{Y}.

In the paper we prove an analog of Theorem \ref{as} for positive
contractions of $L^1(M,\t)$ associated with a finite von Neumann
algebra. We think that this theorem will serve to prove the
existence of a stationary state for given dynamical system. Note
that the existence of a stationary state is an actual problem to
study ergodic properties of the quantum dynamical system (see
\cite{F},\cite{J2}). As a consequence of our main result we infer
that the mixing property of positive $L^1$-contraction implies
completely mixing property one, i.e. we prove a noncommutative
extension of Theorem \ref{ks}.  It should be noted that our
results are not valid when von Neumann algebra is semi-finite.

\section{Preliminary}

Throughout the paper  $M$ would be a von Neumann algebra with the
unit $\e$ and let $\t$ be a faithful  normal finite trace on $M$.
Therefore we omit this condition from the formulation of theorems.
Recall that an element $x\in M$ is called {\it self-adjoint} if
$x=x^*$. The set of all self-adjoint elements is denoted by
$M_{sa}$. A self-adjoint element $p\in M$ is called a {\it
projection} if $p^2=p$. The set of all projections in $M$ we will
denote by $\nabla$. By $M_*$ we denote a pre-dual space to $M$
(see for definitions \cite{BR},\cite{T}).

The map $\|\cdot\|_{1}:M\rightarrow [0, \ \infty)$ defined by the
formula $\|x\|_{1}=\tau(|x|)$ is a norm (see \cite{N}). The
completion of $M$ with respect to the norm $\|\cdot\|_{1}$ is
denoted by $L^{1}(M,\t)$. It is known \cite{N} that the spaces
$L^{1}(M,\tau)$  and $M_*$ are isometrically isomorphic, therefore
they can be identified. Further we will use this fact without
noting.

\begin{thm}\label{2.1}\cite{N} The space $L^{1}(M,\tau)$
coincides with the set
$$
L^{1}=\{ x=\int^{\infty}_{-\infty}\lambda de_{\lambda}\
:\int^{\infty}_{-\infty}|\lambda| d \tau (e_{\lambda})< \infty \}.
$$
Moreover,
$$
\|x\|_{1}=\int^{\infty}_{-\infty}|\lambda |d \tau (e_{\lambda}).
$$
\end{thm}

It is known \cite{N} that the equality $$
L^1(M,\t)=L^1(M_{sa},\t)+iL^1(M_{sa},\t) $$ is valid. Note that
$L^1(M_{sa},\t)$ is a pre-dual to $M_{sa}$.

Let $T:L^1(M,\t)\to L^1(M,\t)$ be linear bounded operator.  We say
that a linear operator $T$ is  {\it positive} is $Tx\geq 0$
whenever $x\geq 0$.  A linear operator $T$ is said to be a {\it
contraction} if $\|T(x)\|_1\leq \|x\|_1$ for all $x\in
L^1(M_{sa},\t)$.

\section{Mixing and completely mixing contractions}

Let $M$ be a von Neumann algebra with faithful normal finite trace
$\tau$. Let $L^1(M,\tau)$ be a $L^1$-space. In the sequel by
$\|\cdot\|$ we mean the norm $\|\cdot\|_1$.

Let $T:L^1(M,\tau)\to L^1(M,\tau)$ be a linear contraction. Define
\be \label{mix}
\bar\rho(T)=\sup\left\{\lim_{n\to\infty}\frac{\|T^n(u-v)\|}{\|u-v\|}
: \ \ u,v\in L^1(M_{sa},\t), u,v\geq 0, \|u\|=\|v\|\right\}. \ee

If $\bar\rho(T)=0$ then $T$ is called {\it completely mixing}.
Note that certain properties  of completely mixing quantum
dynamical systems have been studied in \cite{AP}.

Denote
$$
X=\{x\in L^1(M_{sa},\tau): \ \t(x)=0\}.
$$
Recall that a positive contraction $T$ is {\it mixing} if for all
$x\in X$ and $y\in M$ the following condition holds
$$
\lim_{n\to\infty}\t(T^n(x)y)=0.
$$

Let $T$ be a positive contraction of $L^1(M,\t)$, and let $x\in
L^1(M,\t)$ be such that $x\geq 0$, $x\neq 0$. We say that $T$ is
{\it smoothing} with respect to(w.r.t.) $x$ if for every $\i>0$
there exist $\delta>0$ and $n_0\in\N$ such that $\t(pT^nx)<\i$ for
every $p\in\nabla$ such that $\t(p)<\delta$ and for every $n\geq
n_0$. A commutative analog of this notion was introduced in
\cite{ZZ},\cite{KT}.\\

\begin{thm}\label{as-n} Let $T:L^1(M,\t)\to L^1(M,\t)$ be a
positive contraction. Assume that  there is a positive element
$y\in L^1(M,\t)$ such that $T$ is smoothing w.r.t. $y$. Then
$\lim\limits_{n\to\infty}\|T^ny\|=0$ or there is a non zero
positive $z\in L^1(M,\t)$ such that $Tz=z$.
\end{thm}

\begin{proof} The contractivity of $T$ implies that the limit
$$
\lim_{\n\to\infty}\|T^ny\|=\a
$$
exists. Assume that $\a\neq 0$. Define $\l:M_{sa}\to \R$ by
$$
\l(x)=L((\t(xT^ny)_{n\in\N}))
$$
for every $x\in M_{sa}$, here $L$ means a Banach limit (see,
\cite{K}). We have
$$
\l(\e)=L((\t(T^nx)_{n\in\N}))=\lim_{\n\to\infty}\|T^nx\|=\a\neq 0,
$$
therefore $\l\neq 0$. Besides, $\l$ is a positive functional,
since for positive element $x\in M_{sa}$,$x\geq 0$ we have
$$
\t(xT^ny)=\t(x^{1/2}T^nyx^{1/2})\geq 0,
$$
for every $n\in\N$.

For arbitrary $x\in M$, we have $x=x_1+ix_2$ and define $\l$ by
$$
\l(x)=\l(x_1)+i\l(x_2).
$$

Let $T^{**}$ be the second dual of $T$, i.e. $T^{**}:M^{**}\to
M^{**}$. The functional $\l$ is $T^{**}$-invariant. Indeed,
\bea\label{inv}
(T^{**}\l)(x)=<x,T^{**}\l>=<T^*x,\l>=\nonumber \\
=L((\t(T^nyT^*x)_{n\in\N}))=L((\t(xT^{n+1}y)_{n\in\N}))=\nonumber
\\
=L((\t(xT^ny)_{n\in\N}))=\l(z).\nonumber \eea

Let $\l=\l_n+\l_s$ be the Takesaki's decomposition (see \cite{T})
of $\l$ on normal and singular components. Since  $T$ is normal
and $T^{**}\l=\l$, so using the idea of \cite{J2} it can be proved
the equality $T^{**}\l_n=\l_n$. Now we will show that $\l_n$ is
nonzero. Consider a measure $\m:=\l\upharpoonright_\nabla$. It is
clear that $\m$ is an additive measure on $\nabla$. Let us prove
that it is $\s$-additive. To this end, it is enough to show that
$\m(p_k)\to 0$ whenever $p_{k+1}\leq p_k$ and $p_k\searrow 0$,
$p_k\in\nabla$.

Let $\i>0$. From $p_n\searrow 0$ we infer that $\t(p_n)\to 0$ as
$n\to\infty$. It follows that there exists $k_{\i}\in\N$ such that
$\t(p_k)<\i$ for all $k\geq k_\i$. Since $T$ is smoothing w.r.t.
$y$ we obtain
$$
\t(p_kT^ny)<\i, \ \ \ \ \forall k\geq k_\i,
$$
for every $n\geq n_0$. From a property of Banach limit we get
$$
\l(p_k)=L((\t(p_kT^ny)_{n\in\N})<\i \ \ \ \ \textrm{for every} \ \
k\geq k_\i,
$$
which implies $\m(p_k)\to 0$ as $k\to\infty$. This means that the
restriction of $\l_n$ on $\nabla$ coincides with $\m$. Since
$$
\t(p^{\perp}T^ny)>\t(T^ny)-\i\geq\inf\|T^ny\|-\i=\a-\i
$$
as $\i$ has been arbitrary, so $\a-\i>0$, and hence
$\m(p^{\perp})>0$ for all $p\in\nabla$ such that $\t(p)<\delta$.
Therefore $\m\neq 0$, and consequently, $\l_n\neq 0$.

From this we infer that there exists a positive element $z\in
L^1(M,\t)$ such that
$$
\l_n(x)=\t(zx), \ \ \ \forall x\in M.
$$
The last equality and $T^{**}\l_n=\l_n$ yield that  \bea
\t(zx)=<x,T^{**}\l_n>=<T^*x,\l_n>=\t(zT^*x)=\t(Tzx) \nonumber \eea
for every $x\in M$, which implies that $Tz=z$. \end{proof}

Using pre-compactness criterion for a subset of $L^1(M,\t)$ (see
\cite{T}) one  can prove the following

\begin{lem}\label{wek} Let $(x_n)\subset L^1(M,\t)$,
$\sup\limits_{n\in\N}\|x_n\|<\infty$. Assume that $x_n\to x^*$
weakly. Then for an arbitrary $\i>0$ there exist $\delta>0$ and
$n_0\in\N$ such that $\t(p|x_n|)<\i$ for every $p\in\nabla$ such
that $\t(p)<\delta$ and $\forall n\geq n_0$.
\end{lem}

\begin{cor}\label{cmp} Let $x\in L^1(M,\t)$, $x\geq 0$.
Assume that $T^nx\to x^*$ weakly. Then $T$ is smoothing w.r.t.
$x$.
\end{cor}

{\it Remark 3.1.} The proved Theorem \ref{as-n} is a
non-commutative analog of Akcoglu and Sucheston result \cite{AS}.
But they used weak convergence instead of smoothing. Lemma
\ref{wek} shows that smoothing condition is less restrictive than
the one they used.\\

Before proving next theorem let us give the following an auxiliary

\begin{lem}\label{pos} Let $x\in L^1(M,\t)$. If  the
inequality \be\label{zero} \t(xy)\geq 0\ee is valid for every
$y\geq 0$, $y\in M$. Then $x\geq 0$.
\end{lem}

\begin{proof}. The following is given $x=x^+-x^-$. Let
$$
x=\int_{-\infty}^{\infty}\lambda de_{\lambda}
$$
be the spectral resolution of $x$. Set
$$
p=\int_{-\infty}^{0}de_{\lambda}.
$$

Then according to (\ref{zero}) one gets $\t(xp)\geq 0$. On the
other hand we have $xp=-x^-$, hence $\t(x^-)\leq 0$, since
$x^-\geq 0$ and $\t$ is faithful we infer that $x^-=0$. Therefore
$x=x^+\geq 0$. \end{proof}

From Lemma \ref{wek} and Theorem \ref{as-n} we find the following

\begin{thm}\label{ks-n} Let $T:L^1(M,\t)\to L^1(M,\t)$ be a
positive contraction such that $ |T(x)|\leq T(|x|)$ for every
$x\in L^1(M,\t)$,$x=x^*$.  Assume that  there exits no non zero
$y\in L^1(M,\t)$, $y\geq 0$ such that $Ty=y$. If  for $z\in
L^1(M,\t)$ the sequence $(T^nz)$ converges weakly to some element
of $L^1(M,\t)$, then $\lim\limits_{n\to\infty}\|T^nz\|=0$.  In
particular, if $T$ is mixing, then $T$ is completely mixing.
\end{thm}

\begin{proof} As in the proof of Theorem \ref{as-n} we assume that
$$
\lim_{n\to\infty}\|T^nz\|=\a>0.
$$

Define $\l:M_{sa}\to \R$ by
$$
\l(x)=L((\t(x|T^nz|)_{n\in\N}))
$$
for every $x\in A$. Using the same argument as in the proof of
Theorem \ref{as-n} one can show that there exists nonzero positive
element $y\in L^1(M,\t)$ such that
$$
\l_n(x)=\t(yx), \ \ \ \forall x\in M.
$$
here $\l_n$ is the normal part of $\l$.

From the property of $T$ we infer \bea
\t(Tyx) &= &\t(yT^{*}x)\nonumber \\
& = &L((\t(|T^nz|T^*x))_{n\in\N})\nonumber \\
& = & L((\t(T|T^nz|x))_{n\in\N})\nonumber \\
& \geq & L((\t(|T^{n+1}z|x))_{n\in\N})=\t(yx)\nonumber \eea for
all $x\geq 0$. Hence, for every $x\geq 0$ we have
$$
\t((Ty-y)x)\geq 0.
$$

According to  Lemma \ref{pos} we infer that  $Ty\geq y$. Since $T$
is a contraction one gets $Ty=y$. But this contradicts to the
assumption of theorem. \end{proof}

{\it Remark 3.2.} The proved theorem is an non-commutative analog
of Theorem \ref{ks}. Certain similar results has been obtained in
\cite{L},\cite{FR2} for quantum dynamical semigroups in
von Neumann algebras.\\

\begin{cor}\label{jor}  Let $\a:M\to M$ be a normal Jordan
automorphism such that  there exits no non zero $y\in L^1(M,\t)$,
$y\geq 0$ such that $\a^*y=y$, where $\a^*$ is the conjugate
operator to $\a$.  If for $z\in L^1(M,\t)$ the sequence
$((\a^*)^nz)$ converges weakly to some element of $L^1(M,\t)$,
then $\lim\limits_{n\to\infty}\|(\a^*)^nz\|=0$.
\end{cor}

The proof immediately comes from Theorem \ref{ks-n} since for
Jordan automorphisms the equality $|\a(x)|=\a(|x|)$ is valid for
all $x\in M_{sa}$ (see \cite{BR}).\\

{\it Remark 3.3.} Note that analogous theorem has been recently
proved by A.Katz \cite{Ka} for automorphisms of an arbitrary von
Neumann algebra. Corollary \ref{jor} extends his result to Jordan
automorphisms of finite von Neumann algebras. Here it should be
also noted that linear mappings of von Neumann algebras which
satisfy the condition $|\a(x)|=\a(|x|)$ have been studied in
\cite{Ra}.\\

{\it Remark 3.4.} It should be noted that  Theorems \ref{as-n} and
\ref{ks-n} are not valid if a von Neumann algebra is semi-finite.
Indeed, let $B(\ell_2)$ be the algebra of all linear bounded
operators on Hilbert space $\ell_2$. Let $\{\phi_n\}$,$n\in\N$ be
a standard basis of $\ell_2$, i.e.
$$
\phi_n=(\underbrace{0,\cdots,0,1}_n,0\cdots).
$$
The matrix units of $B(\ell_2)$ can be defined by
$$
e_{ij}(\xi)=(\xi,\phi_i)\phi_j, \ \ \ \xi\in \ell_2, \ i,j\in\N.
$$
A trace on $B(\ell_2)$ is defined by
$$
\t(x)=\sum_{k=1}^\infty(x\phi_k,\phi_k).
$$
By $\ell_{\infty}$ we denote a maximal commutative subalgebra
generated by elements $\{e_{ii}:\  i\in\N\}$. Let $E:
B(\ell_2)\to\ell_{\infty}$ be the canonical conditional
expectation (see \cite{T}). Define a map
$s:\ell_{\infty}\to\ell_{\infty}$ as follows: for every element
$a\in\ell_{\infty}$, $a=\sum\limits_{k=1}^\infty a_ke_{kk}$ put
$$
s(a)=\sum_{k=1}^\infty a_{k}e_{k+1,k+1}.
$$

Define $T:B(\ell_2)\to B(\ell_2)$ as $T(x)=s(E(x))$, $x\in
B(\ell_2)$. It is clear that $T$ is positive and
$\t(T(x))\leq\t(x)$ for every $x\in L^1(B(\ell_2),\t)\cap
B(\ell_2)$, $x\geq 0$. Hence, $T$ is a positive $L^1$-contraction.
But for this $T$ there is no non zero $x$ such that $Tx=x$.
Moreover, for every $y\in L^1(B(\ell_2),\t)$ we
have $\lim\limits_{n\to\infty}\|T^ny\|_1\neq 0$.\\

\section*{acknowledgements}
The first named author (F.M.) thanks TUBITAK-NATO PC-B programme
for providing financial support and Harran University for kind
hospitality and providing all facilities. The authors would like
to thank  Prof. V.I.Chilin from National University of Uzbekistan,
for valuable advice on the subject. The work is also partially
supported by Grant $\Phi$-1.1.2 Rep. Uzb.

The authors also express their gratitude to the referee's helpful
comments.

\bibliographystyle{amsplain}

\end{document}